\documentclass[reqno]{amsart}
\usepackage{graphicx}
\theoremstyle{theorem}
\newtheorem{theorem}{Theorem}[section]
\newtheorem{corollary}[theorem]{Corollary}
\newtheorem{lemma}[theorem]{Lemma}
\newtheorem{prop}[theorem]{Proposition}

\theoremstyle{definition}
\newtheorem{definition}{Definition}[section]
\newtheorem{remark}[definition]{Remark}
\newtheorem{assumption}[definition]{Assumption}
\newtheorem{example}[definition]{Example}

\numberwithin{equation}{section}

\newcommand{\CC}{\mathbb{C}}
\newcommand{\RR}{\mathbb{R}}
\newcommand{\QQ}{\mathbb{Q}}
\newcommand{\ZZ}{\mathbb{Z}}

\newcommand{\FF}{\mathcal{F}}

\newcommand{\CM}{\mathcal{M}}
\newcommand{\MM}{\mathcal{M}}

\newcommand{\CP}{\mathbb{C}P}
\newcommand{\AI}{A_\infty}
\newcommand{\LI}{L_\infty}

\newcommand{\NOV}{C^*(L;\Lambda_{nov})}
\begin{document}
\title[Products of Floer cohomology of toric Fano manifolds]
{Products of Floer cohomology \\
of torus fibers in toric Fano manifolds}

\author{Cheol-Hyun Cho}
\address{Department of Mathematics, Northwestern University,
Evanston, IL 60208, Email: cho@math.northwestern.edu}

\begin{abstract}
We compute the ring structure of Floer cohomology groups of Lagrangian torus fibers in some toric Fano manifolds
continuing the study of \cite{CO}. Related $\AI$-formulas hold for transversal choice of chains.
Two different computations are provided: a direct calculation 
using the classification of holomorphic discs by Oh and the author
in \cite{CO}, and another method by using an {\it analogue of divisor equation} in Gromov-Witten invariants to the
case of discs.
 Floer cohomology rings are shown to be isomorphic to Clifford algebras, whose
quadratic forms are given by the Hessians of functions $W$, which turn out
to be the superpotentials of Landau-Ginzburg mirrors. In the case
of $\CP^n$ and $ \CP^1 \times \CP^1$, this proves the prediction
made by Hori, Kapustin and Li by B-model calculations via physical arguments.
 The latter method also provides correspondence between higher derivatives of
the superpotential of LG mirror with the higher products of $\AI$(or $\LI$)-algebra of
the Lagrangian submanifold. 
\end{abstract}

\maketitle

\bigskip

\section{Introduction}
Floer theory of Lagrangian intersections has been proved to be a
powerful technique in symplectic geometry. Also since the
``homological mirror symmetry'' conjecture by Kontsevich \cite{Ko}, it
has become much more exciting field of mathematics, which yet has a long
way to be fully understood. Recently, Fukaya, Oh, Ohta and Ono constructed 
$\AI$-algebra of Largrangian submanifold and Floer homology in general setting in their beautiful work \cite{FOOO}.
But the construction is highly non-trivial to overcome several technical problems. 
The first problem is the well-definedness of the moduli space of $J$-holomorphic discs compatible for
all homotopy classes. It was observed in \cite{FOOO}, that standard Kuranishi perturbation
does not produce compatible and transversal moduli space in general. Another problem is that
even if moduli spaces of $J$-holomorphic discs are well-defined, it does not directly produce
$\AI$-algebra since one has to work at the chain level. 

In \cite{CO}, Yong-Geun Oh and the author has explicitly described the moduli space
of holomorphic discs in the case of Lagrangian torus fibers in toric Fano manifolds, and used that information to compute Floer cohomology groups. 
A combinatorial discription of a fiber whose
Floer cohomology is non-vanishing was found, and for such a fiber, the Floer cohomology was in fact isomorphic to
singular cohomology as a module. It was shown that
all of holomorphic discs in these cases are transversal. To compactify the moduli space, we need an additional assumption
regarding the behavior of holomorphic spheres on toric Fano manifold(see Assumption \ref{assump}).
In this paper, we first consider related $\AI$-algebra which is defined transversally.
Namely, fiber products with various chains in the Lagrangian submanifold $L$ in the definition of $\AI$-algebra  can be made transversal for the generic choice of chains.
This gives a partial $\AI$-algebra, but products on the cohomology of these $\AI$-algebras are shown to be well-defined. 
How to obtain an actual $\AI$-algebra from this partial algebra is an interesting question. With 
skew-symmetrization in this toric Fano case, these partial $\AI$-algebras gives well-defined $\LI$-algebras.
On the other hand, recently Fukaya has constructed an $\AI$-algebra on DeRham complex of Lagrangian submanifolds. A computation in toric Fano case can be carried out in the DeRham setting, which will produce actual $\AI$-algebra.

 Then we show that Floer cohomology ring $HF^{BM}(L;J_0)$  is isomorphic to a Clifford algebra $Cl(V,Q)$ where
$Q$ is a symmetric bilinear form. It is very interesting that
the symmetric bilinear form $Q$ we obtained
exactly agrees with the Hessian of the superpotential $W$
of the mirror Landau-Ginzburg model studied by Hori and Vafa \cite{HV}. 
(This is related to homological mirror symmetry conjecture between A-model in Fano manifolds and
$B$-model in Landau-Ginzburg mirror.)
In particular, the Floer cohomology of the Clifford torus $T^n$ in $\CP^n$
is isomorphic to the Clifford algebra with $n$ generators as a ring. 

Such product structures in the Clifford torus $T^n$ in $\CP^n$ and $T^1 \times T^1$ in 
$\CP^1 \times \CP^1$ have been
conjectured by Hori and Kapustin and Li \cite{KL}, recently in general by \cite{KL2}
 from the calculation on B-model side using physical arguments.
Mathematical account of the product structure on $B$-model side looks plausible
 considering the paper by Orlov \cite{O}.

We provide two ways of computing the product structure. First, we provide direct computations 
exploiting the classification of all holomorphic discs with boundary on $L$
by Oh and the author (\cite{CO}).
Another method is  by using an analogue of {\it divisor equation} for discs, which is introduced in section
\ref{divisor}. The latter method easily provides the general correspondence
 between higher derivatives of 
the superpotential of LG mirror with the higher products of $\AI$(or $\LI$)-algebra of
Lagrangian submanifold. This extends the correspondence proved by Oh and the author in \cite{CO}
that obstruction cochain $m_0 = l_0$ agrees with the superpotential itself
and non-vanishing of Floer cohomology corresponds to the critical points of the superpotential $W$.
These $l_\infty$-products are invariant under the perturbation of an almost complex structure. 

We also provide an explicit filtered chain map between singular cochain complex and
Bott-Morse Floer complex in the case of torus fibers $L$ in toric Fano manifolds, which induces an isomorphism
in cohomology in case Floer homology is non-vanishing.

{\em Acknowledgements.}
The author do not claim originality of the construction of $\AI$-algebra
of Lagrangian submanifold which should be given to Fukaya, Oh, Ohta and Ono for their 
ingenious work. We would like to thank Yong-Geun Oh for reading
the original draft and for helpful comments. Part of the paper is written during author's visit
to Mathematical Science Research Institute, and he would like to
thank for its hospitality.

\section{$\AI$-algebra of Lagrangian submanifold}
In this section we recall the construction of the $\AI$-algebra of a Lagrangian submanifold. In fact,
we will provide a {\it transversal} version (partial $\AI$-algebra) which is suitable for our purposes.
(This version is only suitable for the case when the moduli space is already well-defined).

The $\AI$-algebra in this case naturally arises from the 
stable map compactification of the moduli spaces of holomorphic discs.
The moduli space of a disc with $n+1$ boundary marked points,
$\CM_{n+1}$, can be seen also as a compactification of a configuration
space of $n-2$ points on an interval $[0,1]$. (By $Aut(D^2)$, send
$n+1,0,1$-st marked points to $1,\infty,0$ where we identify
$D^2$ with the upper-half plane). The latter gives well-known 
Stasheff Polytope \cite{S1}. 

We first recall the definition of (non-unital) $\AI$-algebra introduced
by Stasheff \cite{S1}.
Let $A = \oplus_{i\in \ZZ} A^{i}$ be a $\ZZ$-graded module over $R$,
where $R$ is a commutative ring with unit.
As usual, we denote its suspension by $A[1]^i = A^{i+1}$. 
\begin{definition}
A structure of  (non-unital) $\AI$-algebra on $A$ is given by 
a series of $R$-module homomorphisms $m_k:A^{\otimes n} \to A[2-n]$ for non negative integer $k$,
satisfying quadratic equations
\begin{equation}\label{aiformula}
\sum_{k_1+k_2=k+1} \sum_i (-1)^{deg \; x_1 +\cdots+deg \; x_{i-1} + i-1}
\end{equation}
$$m_{k_1}(x_1,\cdots,m_{k_2}(x_i,\cdots,x_{i+k_2-1}),\cdots,x_k)=0.$$
\end{definition}

In the transversal version, the above formula will only hold 
on a dense transversal sequence of chains for each $k$.

Now we recall the setting for the objects of the chain complex.
We refer readers to \cite{FOOO} Appendix A for a complete explanation about introducing this setup. 
Let $C^*(L;\Lambda_{nov})$ be the set of currents on $L$ realized by
geometric chains as follows:
For a given (n-k)-dimensional geometric chain $[P,f]$, we consider the current
$T([P,f])$ which is defined as follows: The current $T([P,f])$ is
an element in $D'^k(M;\RR)$ where $D'^k(M;\RR)$ is the set of
distribution valued k-forms on $M$ : For any smooth (n-k)-form
$\omega$, we put
\begin{equation}\label{current}
\int_M T([P,f]) \wedge \omega = \int_P f^* \omega
\end{equation}
This defines a homomorphism $$ T: S_{n-k}(M;\QQ) \to D'^k(M;\RR)$$
where $S_{n-k}(M;\QQ)$ is the set of all (n-k) dimensional
geometric chains with $\QQ$-coefficient. Let
$\overline{S}^k(M,\QQ)$ be the image of the homomorphism $T$. 
We extend the coefficient ring $\QQ$ to $\Lambda_{nov}$.
Then we set
\begin{equation}\label{ch}
C^k(L;\Lambda_{nov}) :=  \overline{S}^k(M,\Lambda_{nov})
\end{equation}
Since we consider
the elements in the image of $T$, if the image of the map $f$ of
the geometric chain $[P,f]$ is smaller than expected dimension,
then it gives 0 as a current. This fact will be used crucially
later on. Also, note that the map $T$ is not injective, hence
some elements get identified under the map $T$.
Also note that we take the whole image of $T$ (instead of 
taking a countable subset of it) as transversality
of fiber products in the definition of $m_k$ is achieved by choosing generic chains.

The classical part of the maps $\{m_k\}$ are defined as follows, which
is different from that of \cite{FOOO} (In \cite{FOOO}, $m_{k,0}$ defines
an $\AI$-algebra of singular cochains).
\begin{definition}\label{mk0}
The maps  $m_{k,0}$ for $k=0,1,\cdots$ on $C^*(L;\Lambda_{nov})$
 are transversally defined by the following maps. For $[P,f],[Q,g] \in C^*(L;\QQ)$,
\begin{enumerate}
\item $m_{0,0} =0 $.
\item $m_{1,0} ([P,f]) = (-1)^n [\partial P,f].$
\item $m_{2,0} ([P,f],[Q,g]) = (-1)^{deg P (deg Q +1)}[f(P) \cap g(Q), i]=0.$
where $i$ is an embedding into $L$.
\item for $k \geq 3$, 
\begin{equation} m_{k,0} \equiv 0 \end{equation}
\end{enumerate}
We extend above maps linearly over $\Lambda_{nov}$.
The notation $\partial$ here is the usual boundary operator for singular homology.
\end{definition}

Now, the quantum contribution part is defined in the same way as in \cite{FOOO}.
\begin{definition}\cite{FOOO}\label{mdef}
\begin{enumerate}
\item For a geometric chain $[P,f]\in C^{g}(L:\QQ)$ and non-zero $\beta$, define 
\begin{equation}\label{m0}
m_{0,\beta} = [\CM_1(\beta),ev_0]. 
\end{equation}
\begin{equation}\label{m1}
m_{1,\beta}[P,f] = [\CM_2(\beta)\;_{ev_1} \times _f P, ev_0] 
\end{equation}

\item For each $k \geq 2$, non-zero $\beta$,
for geometric chains $$[P_1,f_1]\in C^{g_1}(L:\QQ), \cdots,
[P_k,f_k]\in C^{g_k}(L:\QQ)$$
(i.e. dimension of $[P_i,f_i]$ as a chain is
$n-g_i$), define 
$$m_{k,\beta}([P_1,f_1],\cdots,[P_k,f_k])$$
\begin{equation}
=(-1)^{\epsilon}[\CM_{k+1}^{\textrm{main}}(\beta) \,_{(ev_1,\cdots,ev_k)}
\times_{(f_1,\cdots,f_2)} (P_1 \times \cdots \times P_k), ev_0]
\end{equation}
Here $\epsilon$ is a sign assigned as follows:
\begin{equation}
\epsilon = (n+1) \sum_{j=1}^{k-1} \sum_{i=1}^j \textrm{deg}(P_i)
\end{equation}
\item Then we define the maps $m_k \;(k \geq 0)$ by
$$m_k([P_1,f_1],\cdots,[P_k,f_k]) = \sum_{\beta \in \pi_2(M,L)} m_{k,\beta}(
[P_1,f_1],\cdots,[P_k,f_k])\otimes
T^{Area(\beta)} q^{\mu(\beta)/2}.$$
\end{enumerate}
\end{definition}
\begin{remark}
Here $\CM_k(\beta)$ is a compactified moduli space of $J$-holomorphic discs
with $k$ marked point on $\partial D^2$.
Recall that $\CM_k(\beta)$ for $k\geq 3$ has several connected component
by the ordering of the $k$ marked points on $\partial D^2$
and by $\CM^{main}_k(\beta)$ we denote the connected component
where marked points $z_1,\cdots, z_k$ lie cyclically on $\partial D^2$ 
counter-clockwise. 
\end{remark}

Also, the fiber products defined above are not always transversal, and we discuss this issue in 
section \ref{sec:trans}

Here we recall the dimension formula of $m_{k,\beta}$ 
when the involved fiber product is transversal.
\begin{prop}[\cite{FOOO} Proposition 13.16]\label{dim} For non-zero $\beta$, when transversal,
$$m_{k,\beta}((P_1,f_1),\cdots,(P_k,f_k)) \in C_{n- \sum_{i=1}^k g_i +
\mu(\beta) -2 +k } (L: \QQ)$$
\end{prop}

\begin{prop}[cf. \cite{FOOO}] These $\{ m_k \}$ maps satisfy the $A_{\infty}$ formulas (\ref{aiformula})
for transversal sequence of of chains in $C^*(L;\Lambda_{nov})$. 
\end{prop}
\begin{proof}
This is essentially the theorem proved in \cite{FOOO}. We recall its proof for the convenience of
readers and explain the changes made for $m_{k,0}$. 

For simplicity, we recall the proof the only the third $\AI$-formula.
Consider the moduli space of J-holomorpic discs intersecting chains $P$ and $Q$
(See Figure \ref{fig1}).
\begin{equation}\label{big}
m_{2,\beta}(P,Q)= (\CM^{main}_3(\beta)_{ev_1,ev_2} \times_{f,g} (P \times Q)
,ev_0)
\end{equation}

Now, we consider all possible stable map compactification of this moduli space
and its image under the evaluation map.
The limit configurations of codimension 1 of the image  can be written as follows.
See Figure \ref{fig1}, where each figure corresponds to the following terms.

\begin{figure}
\begin{center}
\includegraphics[height=2.5in]{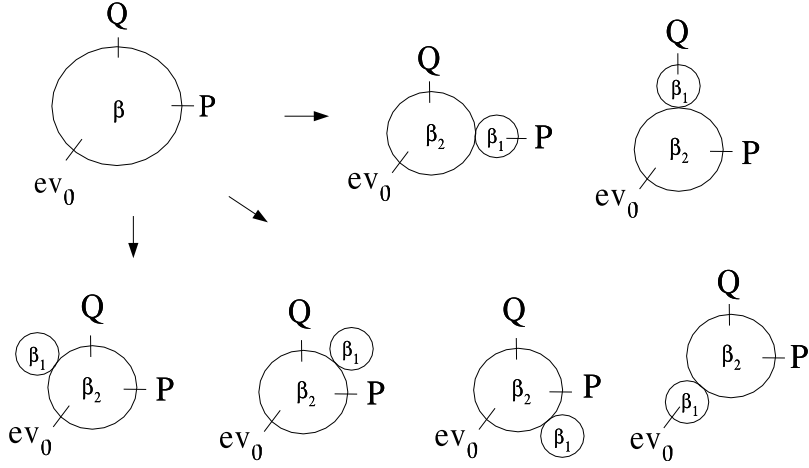}
\caption{Limit configurations of (\ref{big}) of codimension 1.}
\label{fig1}
\end{center}
\end{figure}

\begin{equation}
 m_{2,\beta}(P,Q) \rightarrow m_{2,\beta_2}(m_{1,\beta_1}(P),Q) ,
m_{2,\beta_2}(P,(m_{1,\beta_1}(Q),)
\end{equation}
\begin{equation}
m_{3,\beta_2}(P,Q,m_{0,\beta_1}), 
 m_{3,\beta_2}(P, m_{0,\beta_1},Q), 
m_{3,\beta_2}(m_{0,\beta_1},P,Q), 
m_{1,\beta_1}(m_{2,\beta_2}(P,Q))
\end{equation}

Degenerations into several (three or more) disc components or
sphere bubbles also occur.
But if transversalities are satisfied for such singular strata
with positivity assumptions on Lagrangain submanifold,
such strata should be of codimension 2 or more, hence they do not contribute
to the $\AI$ formulas.

Now, these limit configurations can be written into an $\AI$-formula up to sign:
$$ \partial (m_{2,\beta}(P,Q)) =  \pm
m_{2,\beta_2}(m_{1,\beta_1}(P),Q) \pm
m_{2,\beta_2}(P,(m_{1,\beta_1}(Q),)$$
$$\pm m_{3,\beta_2}(P,Q,m_{0,\beta_1})
\pm m_{3,\beta_2}(P,m_{0,\beta_1},Q) 
\pm m_{3,\beta_2}(m_{0,\beta_1},P,Q) 
\pm m_{1,\beta_1}(m_{2,\beta_2}(P,Q)) $$

This is the third $\AI$-formula in (\ref{aiformula}) up to sign, and
other formulas can be obtained in a similar fashion by
choosing a $m_{k,\beta}(P_1,\cdots,P_k)$ for general $k$ in (\ref{big}). 

Now we justify the changes made in the definitions of $m_{k,0} \equiv 0$ for $k \geq 3$. 
Consider one of the $m_{3,0}$ term appeared  in the above configuration
when $\beta_2 =0$, or generally one may consider the geometric chain 
$m_{3,0} (P,Q,R)$.
The dimension of the image under the evaluation map of
 $m_{3,0}(P,Q,R)$ is always smaller than the virtual dimension of the
 moduli space: The reason is that
evaluation map of a constant disc forgets the moduli parameter. 
Namely, before evaluation, there is a parameter describing the position of
four marked points on a disc. Recall that 
moduli space of 4 marked points on $\partial D^2$ up to automorphisms of $D^2$ is diffeomorphic to $\RR $ (See \cite{FO}).
But as we evaluate on a constant disc, the image is always a point,
while the moduli parameter is lost under the evaluation map.
Hence, such a term $m_{3,0}(P,Q,R)$ is of codimension 1 
by the virtual dimension, but its actual image is of codimension 2.
Hence terms involving $m_{3,0}$ do not appear in the $\AI$ formula,
which is obtained by considering the codimension 1 boundary
of the image of the chain (\ref{big}) under evaluation map.

This phenomenon always happens for $m_{k,0}$ for any $k \geq 3$
becuase of the same reason. Hence we may set (transversally)
$$ m_{k,0} \equiv 0 \; \; \textrm{for} \; k \geq 3. $$ 
Note that in \cite{FOOO}, the evaluation maps of constant homotopy class 
are also perturbed by moduli parameters, so that the image has the same dimension as virtual dimension
unlike our setting. Also note that $m_{2,0}, m_{1,0}$ does not vanish as there are no moduli parameters in these cases.
This proves the proposition.
\end{proof}

Now, because of the presence of $m_0$ terms, $m_1^2 =0$ does not always hold.
Hence, Floer homology groups are not well-defined in general.
Obstructions for the well-definedness of Floer cohomology was studied in
\cite{FOOO}. In an unobstructed case, one can deform the chain complex
suitable way so that $m_1^2 =0$ holds. For the case of torus fibers in toric Fano manifold,
it is (weakly) obstructed, in which case Floer cohomology with itself is well-defined. 

Related phenomenan in the language of $\AI$-algebra is that 
$m_0$ terms disappear from the $\AI$-formula. 

\begin{prop}[compare,\cite{FOOO} Proposition 7.1]\label{unit}
Let $x_i$ be an element in $C^*(L;\Lambda_{nov})$ for $i=0,\cdots,k$,
for a Lagrangian torus fiber $L$ in toric Fano manifolds. Then, when transversal, we have
$$m_{k+1}(x_1,\cdots,[L],\cdots,x_k)=0, k \geq 2, k=0.$$
$$m_2([L],x_0) = (-1)^{deg (x_0)} m_2(x_0,[L]) = x_0.$$
Namely, [L] behaves as a strict unit.
\end{prop}
\begin{proof}
This was proved in \cite{FOOO}, except that we do not need to
use homotopy unit argument. Recall that  
the {\it forget} maps  commutes with the evaluation maps obviously
in our case, whereas they do not commute in \cite{FOOO} because
of the perturbation of evaluation maps at marked points.
We give the proof of the proposition here for the convenience of readers.
Note that proposition holds by the definition of $m_{k+1,0}$ for $k\geq 2$.
Hence it is enough to show that 
$m_{k,\beta}(x_1,\cdots,[L],\cdots,x_k)=0 $ for $k\geq2$ 
with non-zero $\beta \in \pi_2(M,L)$, and 
the statement about $m_{2,0}$.

Note that the condition for the image of a marked point to meet the fundermental chain $[L]$ is redundant since it always meets $L$.
Hence,
\begin{equation}
ev_0\big(\CM_{k+2}(\beta)\,_{ev} \times (x_1 \times \cdots \times [L] \times \cdots \times x_k) \big)
\subset ev_0( \CM_{k+1}(\beta) \,_{ev}\times (x_1 \times \cdots \times x_k)).
\end{equation}
The dimension of RHS is 
$$ n- \sum deg(x_i) + \mu(\beta)-2 +k.$$
where as the virtual(expected) dimension of LHS is 
$$ n - \sum deg(x_i) +  \mu(\beta) -2 + k+1. $$
Hence, actual image has smaller dimension than expected dimension,
which becomes zero in the language of currents.
The case of $m_{2,0}$ follows from sign convention of \cite{FOOO}.
\end{proof}

The second formula implies that the Floer cohomology is well-defined in this case
as observed in \cite{CO} and \cite{C} Proposition 3.18;
In this case it was shown that $m_0(1)= \sum_{i=1}^N [L] \otimes T^{e_i} q$ is a multiple
of fundermental chain. 

Here, we write the first three $\AI$-formulas where $m_0$ terms are
dropped because of the above proposition.
\begin{eqnarray}\label{AIformula}
0 &=& m_1 \circ m_1  \\
0 &=& m_2(m_1(x),y) + (-1)^{deg(x) +1} m_2(x,m_1(y)) + m_1(m_2(x,y))  \\
0 &=& m_1(m_3(x,y,z)) + m_2(m_2(x,y),z) +(-1)^{deg(x) +1} m_2(x,m_2(y,z))\\
& & \nonumber + m_3(m_1(x),y,z)+ m_3(x,m_1(y),z) +  m_3(x,y,m_1(z))
\end{eqnarray}
The first equation implies that $m_1$ defines the cochain complex. The second equation
implies that $m_2$ defines a product of the cohomology up to sign.
For $x,y,z \in HF^{BM}(L;J_0)$, we have $m_1(x)=m_1(y)=m_1(z)=0$.
Therefore the third equation implies the associativity of the
product up to sign.
\begin{eqnarray}\label{aalgebra}
m_2(m_2(x,y),z) +(-1)^{deg(x) +1} m_2(x,m_2(y,z)) =0
\end{eqnarray}

To define an associative product (with correct sign) on cohomology, one should make
the following change of signs.
\begin{definition}\label{sign}
We define 
\begin{eqnarray}
\widetilde{m}_1(P) &=& (-1)^{degP}m_1(P) \\
\widetilde{m}_2(P,Q) &=& (-1)^{deg P (deg Q +1)} m_2(P,Q).
\end{eqnarray}
\end{definition}
\begin{remark}
The first sign appears due to cohomological sign convention.
The second sign appears due to the sign convention
of \cite{FOOO}
\end{remark}

The resulting $\AI$-formulas for the new $\{ \widetilde{m}_k \}$ are
\begin{eqnarray*}
\widetilde{m}_1(\widetilde{m}_2(x,y)) &=& 
\widetilde{m}_2(\widetilde{m}_1(x),y) + (-1)^{deg \;x}\widetilde{m}_2(x,\widetilde{m}_1(y)) \\
\widetilde{m}_2(\widetilde{m}_2(x,y),z) &=& \widetilde{m}_2(x,\widetilde{m}_2(y,z))
\end{eqnarray*}
for $x,y,z \in HF^{BM}(L;J_0)$. 
Hence $\widetilde{m_2}$ defines an graded associative product on $HF^{BM}(L;J_0)$.

For example, with the new sign, the classical cup product part of
$\widetilde{m}_2$  can be written as
$$\widetilde{m}_{2,0} (P_1, P_2) = P_1 \cap P_2.$$
Also associativity in the classical level is just
$$(P_1 \cap P_2)\cap P_3 =P_1 \cap (P_2\cap P_3).$$

\section{Transversality}\label{sec:trans}
In this section, we discuss the issues regarding the moduli space of $J$-holomorphic discs and
the transversality of $\AI$-algebra.
\subsection{Moduli spaces}
We first recall the following theorem.
\begin{theorem}[\cite{CO}]\label{regularity}
Holomorphic discs in toric manifolds with boundary on any Lagrangian torus fiber are Fredholm regular, i.e.,
its linearization map is surjective.
\end{theorem}
Hence the moduli space of holomorphic discs (before compactification) is a manifold of expected dimensions.
As we try to compactify the moduli space, we may have strata with sphere bubbles. In general toric Fano manifolds,
it is already known that holomorphic spheres are not always Fredholm regular. Hence 
in the compactification of holomorphic discs, some strata (with sphere bubble) may not have the
expected dimension. But since we only evaluate only at the boundary of the discs (not on spheres),
with Fano condition, the evaluation image of such strata is always of codimension of two or higher.
Hence, it is plausible that these moduli spaces with evaluation maps define  currents on $L$. But to make this precise seems
to be a non-trivial problem. Similar problem also has been observed in the case of Gromov-Witten theory
if one try to integrate forms over pseudo-cycle(See page 277 of \cite{MS}). The author do not know how to prove it, so we  require the following strict assumption on sympletic manifold so that the moduli
chain defines a current.
\begin{assumption}\label{assump}
The toric Fano manifold $M$ is assumed to be convex. Namely we require that for any
genus 0 stable map $f:\Sigma \to M$, $f^*T_M$ is generated by global sections.
\end{assumption}
Such assumption holds in the case of complex projective spaces, and products of complex 
projective spaces. Except this rectifiability problem
of the compactified moduli chain of holomorphic discs, the results in this paper holds for all toric Fano manifolds.
Even when the assumption is not satisfied, the results in section \ref{divisor} can be understood independently 
as computations of some invariants. (See Proposition \ref{invariant}).

We remark about perturbing standard complex structure to a tame almost complex structure.
McDuff and Salamon \cite{MS} showed that for a subset $J_{reg}(M)$ of second category, 
so that the moduli spaces of simple $J$-holomorphic curves become pseudo-cycles.
In the case of J-holomorphic discs, it is more complicated since the structure of non-simple J-holomorphic discs can
be very complex. But due to the structure theorem proven by Kwon and Oh \cite{KO}, the similar proof
as in \cite{MS} can be used to show that the moduli space of simple discs are ``pseudo-chain'' which may be
similarly defined as pseudo-cycle. But
also in this case, we do not know if these moduli chains would define currents. If these define currents,
one can prove the invariance of Floer cohomology ring in a similar way as in \cite{FOOO}.

Another approach would be to consider Kuranishi structure of the moduli space of $J$-holomorphic
discs (\cite{FOOO},\cite{FOno}). But as pointed out in \cite{FOOO}, it is not (yet) possible 
to find a Kuranishi perturbation which is compatible for all homotopy classes in $\pi_2(M,L)$. Such
compatibility is rather essential since we are interested in the relations between moduli spaces
which produce $\AI$-formula. 

\subsection{Transversal $\AI$-algebra}
Now we explain how to achieve transversality of the fiber product in the
definition of $\AI$-formulas.
First, recall that ordinary intersection product in the chain level is not well-defined, while cup
product is well-defined on cohomology. Hence, even in the classical level, $\AI$-algebra $(\NOV,m_{k,0})$
is not easy to define, since operations are defined in the chain level. But it is obvious how
to define it to work only transversally. Similar problem occurs for $m_{k,\beta}$. For example
the fiber product $m_k(P,P,\cdots,P)$ is not transversal if $P \neq L$. Hence, authors of
\cite{FOOO} develop non-trivial technique to overcome such problem. In this section,
we show that if we choose the generic sequence of chains, then the fiber product is transversal,
and this transversal $\AI$-algebra  is enough to determine homology and its ring structure. 

\begin{definition}
A $k$-tuple $(P_1,\cdots, P_k)$ is called a {\it transversal sequence}
if the chain $(P_1 \times \cdots \times P_k)$ is transversal
to the image of the map $ev_\beta$ for all $\beta \in \pi_2(M,L)$.
For a transversal sequence $(P_1,\cdots, P_k)$, the fiber product
$m_k(P_1,\cdots, P_k)$ is well-defined.
\end{definition}

Recall that a {\it residual} subset of a space $X$ is one which contains
the intersection of countably many dense open subsets.

\begin{lemma}
For a residual set of $\NOV \times \cdots \times \NOV$,
$k$-th $\AI$-formula (\ref{aiformula}) is well-defined.
Namely the all the fiber products given in the formula are transversal.
\end{lemma}
\begin{proof}
It is enough to show that transversality of 
the chain $(P_1 \times \cdots \times P_k)$ and
the image of $ev_{\beta}$ from each codimension 1 strata of the moduli space
of $J$-holomorphic discs for all $\beta \in \pi_2(M,L)$, which can be 
achieved by choosing generic chains $P_i$'s by the standard transversality theorem.
\end{proof}

\begin{corollary}
 $(\NOV, \{ m_k \})$ satisfies $\AI$-formula for
dense transversal sequence of chains.
\end{corollary}

In fact, in our case it is easy to perturb $(P_1,\cdots,P_k)$ to a transversal sequence due to the presence of torus action.
Namely, as torus $(S^1)^n$ acts on $L$ transitively,
Hence, for a generic $(t_1,\cdots,t_k) \in (S^1)^n \times \cdots \times (S^1)^n$,
$(t_1 \cdot P_1) \times \cdots \times (t_k \cdot P_k) $ is
a transversal sequence. 
Also because we have the same torus action on the moduli space of
holomorphic discs, we have the following identity.
\begin{equation}\label{taction}
m_k (t\cdot P_1,\cdots,t \cdot P_k) = t \cdot m_k(
P_1,\cdots,P_k).
\end{equation}
Therefore, the transversality of $\AI$-formula also can be achieved
by the torus action on each chains: If $m_{k_2}$ term causes non-transversality 
to define $m_{k_1}$ 
in the $\AI$-formula, then we can perturb all chains inside $m_{k_2}$ by
the same $t \in (S^1)^n$ to make $m_{k_2}$ term transversal in $m_{k_1}$ by 
the equality (\ref{taction}). 
Also, it is easy to perturb a Floer-cycle in its cohomology class
by the following lemma.
\begin{lemma}
Let $\Psi$ be the chain map constructed in Definition \ref{correction}.
For any cycle $P$ of singular homology, $\Psi(P)$ is a Floer-cycle.
i.e. $m_1(\Psi(P)) =0$. 
Then for $t \in T^n$, $ t \cdot \Psi(P)$ is
also a Floer cycle, and we have
$$ \Psi(P) - t \cdot \Psi(P) = (-1)^n m_1 \Psi (H)$$
where homotopy $H$ is a singular chain with $m_{1,0} (H) = P - t \cdot P$.
\end{lemma}
\begin{proof}
The equation (\ref{taction}) for $k=1$ implies that
$$m_{1,\beta} (t\cdot P) = t \cdot m_{1,\beta} P.$$
Hence the theorem follows. The last statement follows by applying the Proposition
\ref{iso} (1) for the chain $H$ with the fact that $t \cdot \Psi(P) = \Psi(t \cdot P)$.
\end{proof}

\begin{prop}
$\widetilde{m_2}$ defines a product on the Floer cohomology ring $HF^{BM}(L;J_0)$.
\end{prop}
\begin{proof}
To show that the product is well-defined on cohomology, it is enough to
show that for $P,Q \in \NOV$ with $m_1(P) = m_1(Q) = 0$,
we have $$ m_2(P,t_1 \cdot Q) = m_2(P,t_2\cdot Q) + m_1(R).$$
for generic $t_1,t_2 \in (S^1)^n$ and for some $R \in \NOV$.
First, note that the the fiber product in $m_1(P)$ of $\AI$-algebra is 
transversal for any chain $P$ since the evaluation map from the moduli space is always submersive due to
the torus action. And $m_2(P,Q)$ is transversal if $m_1(P)$ is transversal to $Q$.
Then, for a generic $t \in (S^1)^n$, $m_1(P)$ is transversal to $ t\cdot Q$.
Also, for generic $t_1,t_2 \in (S^1)^n$, $m_1(P)$ is transversal to $H$ with
$m_1(H) = t_1 \cdot Q - t_2 \cdot Q$. (If not, we can perturb $t_1 \cdot Q, t_2 \cdot Q,H$ by
another $t \in (S^1)^n$ to make them transversal.)
Therefore, 
$$m_2(P,t_1 \cdot Q) - m_2(P,t_2 \cdot Q) = m_2(P,m_1(H)) = \pm m_1(m_2(P, H)$$
This finishes the proof.
\end{proof}

\section{Bott-Morse Floer cycles}
In \cite{C} and \cite{CO}, Oh and the present author have shown that
for any such torus fiber $L \subset M$, the Floer homology group $HF(L,L)$
when nonvanishing, is  isomorphic to the singular cohomology of the Lagrangian
submanifold $H^*(L:\Lambda_{nov})$.
Now, we fix a Lagrangian torus fiber $L$ whose Floer cohomology is non-vanishing.
The fact that $HF^{BM}(L;J_0)$ and $H^*(L;\Lambda_{nov})$ is isomorphic as a module is a little bit deceiving
because a cycle in the singular homology is {\it not} a cycle in Floer homology.
We need to modify a cycle, say $P$, by adding correction terms, say $Q$
 to make it satisfy $m_1 (P + Q) =0$.
In the computations of \cite{C} or \cite{CO}, it was automatically taken care of by
the spectral sequence.
We will find exact correction terms for any cycle in proposition \ref{iso}.
Actually we will construct a filtered chain map from singular chain complex to
Bott-Morse Floer complex.

We start with the following definition and an important example
to understand the construction that follows.
\begin{definition}
An element $P = \sum_{i=1}^k a_i\;[P_i,f_i]\; T^{e_i} \;q^{\mu_i} \in \NOV$
is called a Floer-cycle if $m_1(P) =0$.
\end{definition}

\begin{example}\label{pt}
Consider a Clifford torus $T^2$ in $\CP^2$.
A point $<pt>$ is a cycle in  the singuler homology of $T^2$.
Let $l_0,l_1,l_2$ be the cycles in $T^2$ which
are boundaries of holomorphic discs $[z;1;1],[1;z;1],[1;1;z]$.
These three discs have the same symplectic area which we denote
by $\omega(D)$.

Recall from \cite{C} that we have
$$m_1 <pt> = (-1)^n (l_0 + l_1 + l_2)\otimes T^{\omega(D)}q \neq 0.$$
Therefore $<pt>$ is not a Floer-cycle. But, 
$l_0 + l_1 + l_2$ is  homologous to zero. 
We may choose a 2-chain $Q\subset L$ with $ \partial Q =-( l_0+l_1+l_2)$.
Hence $<pt> + Q\otimes T^{\omega(D)}q $ turns out to be a correct Floer-cycle:
\begin{equation}
m_1 (<pt> + Q\otimes T^{\omega(D)}q ) = m_{1,2} (<pt>) +  m_{1,0} (Q)
\otimes T^{\omega(D)}q
\end{equation}
$$ = (-1)^n((l_0 + l_1 +l_2 ) +  \partial Q )\otimes T^{\omega(D)}q = 0.$$
\end{example}

Similarly, we can explicitly construct correction terms as follows
for the general toric Fano case.
We first recall the usual product structure
on the torus $T^n = (S^1)^n$.
$i.e. \; \;$ for $(a_1,\cdots,a_n) \in T^n$,
$(b_1,\cdots,b_n) \in T^n$, we have
$$ (a_1,\cdots,a_n) \times (b_1,\cdots,b_n) = (a_1b_1,\cdots,a_nb_n).$$
Also for  subsets $P \subset T^n, Q \subset T^n$,
we denote by $P \times Q$
$$P \times Q : = \{ (p\times q ) \in T^n | p \in P, q \in Q \}.$$
We may assign the set $P \times Q$ a product orientation.

Recall from \cite{CO} that we have $N$ holomorphic discs of
Maslov index 2 (up to $Aut(D^2)$) with boundary on the Lagrangian
torus fiber $L\subset M$, which we denote by
$D_1, \cdots, D_N$. We denote the homotopy classes of
such discs as $\beta_1, \cdots, \beta_N$.
Then we have 
\begin{equation}\label{mm}
m_{1,\beta_i}(P) = (-1)^n (\partial D_i) \times P. 
\end{equation}
Now, we recall the partition
$$\{1,2,\cdots,N\} = \coprod_{i=1}^l I_i .$$
with respect to the symplectic energy of discs.
i.e. discs $D_j$ for $j \in I_i$ have the same symplectic
area, which we denote as $e_i$.
Nonvanishing of Floer cohomology was shown to be equivalent to
the following equality for each $i=1,\cdots,l$.
$$ \left[ \sum_{j \in I_i} \partial D_j \right] = 0 \;\; \textrm{ in }\; H^*(T^n) $$
\begin{definition}\label{Qi}
For each $i$, we denote by $Q_i$ a 2-chain with the following property.
\begin{equation}\label{qq}
\partial Q_i = - \sum_{j \in I_i} \partial D_j
\end{equation}
We may choose such a 2-chain since RHS is homologus to zero.
\end{definition}

Now, consider the chain complex $\NOV$ defined in (\ref{ch}) with two
different coboundary operators $m_{1,0}$ and $m_1$.
To distinguish two chain complex, we label them as 
$(C_1^*(L,\Lambda_{nov}),m_{1,0})$, whose cohomology is isomorphic to singular cohomology,
and $(C_2^*(L,\Lambda_{nov}), m_1)$, whose cohomology is Bott-Morse Floer cohomology.
Now we define a chain map between these two complexes when
Floer cohomology is non-vanishing.
\begin{definition}\label{correction}
Let $P \subset L$ be any singular chain.
Define
$$\Psi (P) := P + \sum_{i=1}^l (Q_i \times P) \otimes T^{e_i}
+ \sum_{i < j} (Q_i \times Q_j \times P) \otimes T^{e_i + e_j}q^2
+ \cdots $$
$$
+ \sum_{ i_1 < \cdots < i_{k}}
(Q_{i_1} \times \cdots \times Q_{i_k}\times P) \otimes
T^{\sum_{j=1}^k e_{i_j}}q^k + \cdots
+ (Q_1 \times Q_2 \times \cdots \times Q_l \times P) \otimes
T^{\sum_{i=1}^l e_i} q^l .$$
By extending linearly over $\NOV$, we obtain a map
$$\Psi : C_1^*(T^n;\Lambda_{nov}) \to
C_2^*(T^n;\Lambda_{nov}).$$
\end{definition}
\begin{remark}
For simplicity, we define $\Psi$ for singular chains rather than 
geometric chains. It can be easily modified to the latter case.
We also recall that $\NOV$ has a filtration with respect
to energy:
$$\FF^{\lambda_0} C^* = \{\sum_i a_i[P_i,f_i]T^{\lambda_i}q^{m_i} |\lambda_i \geq \lambda_0 \;\textrm{for all}\; i \}$$ 
\end{remark}

\begin{prop}\label{iso}
Let $L$ be a Lagrangian torus fiber in toric Fano manifolds, whose
Floer cohomology is non-vanishing.
Then, the map $\Psi$ defines a filtered chain map
which induces an isomorphism on cohomology.
 $$\Psi : H^*(L;\Lambda_{nov}) \to HF^{BM}(L;J_0).$$
More precisely, 
\begin{enumerate}
\item $ \Psi(m_{1,0} P) = m_1 \Psi(P)$
\item $\Psi (\FF^\lambda (C_1)) \subseteq \FF^\lambda (C_2)$
\end{enumerate}
\end{prop}
\begin{remark}
Note that $\Psi$ is only defined when Floer homology is non-vanishing
since otherwise we can not find chains $Q_i$ in \ref{Qi}.
\end{remark}
\begin{proof}
The second property is clear from the definition, hence we only prove the
first statement, which we prove by direct calculation.
Recall that $m_{1,k} \equiv 0$ for $k\geq 4$ in toric Fano case (see
Proposition 7.2 of \cite{CO}).
Hence, 
\begin{equation}\label{twoeq}
m_1(\Psi(P)) = m_{1,0} \Psi(P) + m_{1,2} \Psi(P).
\end{equation}
The first component can be written as 
\begin{eqnarray*}
m_{1,0} \Psi(P) &=& (-1)^n \partial \Psi(P) \\
&=& (-1)^n(\partial P + \sum_{i=1}^l \partial(Q_i \times P) \otimes T^{e_i} + \cdots ) \\
&=& (-1)^n(\Psi(\partial P) + \sum_{i=1}^l \partial(Q_i) \times P \otimes T^{e_i} +\\
&&
\sum_{i < j} (\partial (Q_i \times Q_j) \times P) \otimes T^{e_i + e_j}q^2 + \cdots).
\end{eqnarray*}
We used the following formula in the last equality, where there is
no sign contribution since $Q_i$'s are 2-chains:
$$\partial (Q_{i_1} \times \cdots \times Q_{i_k}\times P) =
(Q_{i_1} \times \cdots \times Q_{i_k})\times \partial P + 
\sum_{j=1}^{k} (Q_{i_1} \times \cdots (\partial Q_{i_j}) \times Q_{i_k})
\times P. $$

For the second component in (\ref{twoeq}),

\begin{eqnarray*}
 m_{1,2}
 &\sum&_{i_1 < \cdots < i_{k-1}} (Q_{i_1} \times \cdots \times Q_{i_{k-1}}
\times P) \otimes T^{\sum_{l=1}^{k-1} e_{i_l}}q^{k-1} \\
& =&  \sum_{j=1}^N m_{1,\beta_j}
(\sum_{i_1 < \cdots < i_{k-1}} (Q_{i_1} \times \cdots \times Q_{i_{k-1}}
\times P)) \otimes T^{e_j} T^{\sum_{l=1}^{k-1} e_{i_l}}q^k \\
&=&  \sum_{i=1}^l (-(-1)^n\partial Q_i) \times
(\sum_{i_1 < \cdots < i_{k-1}} (Q_{i_1} \times \cdots \times Q_{i_{k-1}}
\times P)) \otimes T^{e_j} T^{\sum_{l=1}^{k-1} e_{i_l}}q^k \\
&=& - (-1)^n\sum_{i_1 < \cdots < i_{k} }
\sum_{j=1}^{k} (Q_{i_1} \times \cdots (\partial Q_{i_j}) \times Q_{i_k}
\times P) \otimes T^{\sum_{l=1}^k e_{i_l}}q^k\\
&=& \sum_{i_1 < \cdots < i_{k} }
 (-(-1)^n)\partial (Q_{i_1} \times \cdots \times Q_{i_k})\times P
\otimes T^{\sum_{l=1}^k e_{i_l}}q^k.
\end{eqnarray*}
In the third equality, we used the identity (\ref{mm}),(\ref{qq})
Hence, we have 
$$m_1(\Psi(P)) = m_{1,0} \Psi(P) + m_{1,2} \Psi(P) = (-1)^n \Psi(\partial P)
= \Psi(m_{1,0}P).$$
\end{proof}

The arguments in this section (hence of the whole paper)
can be extened to the case with different spin structures.
Extension to the case with flat bundles over Lagrangian submanifold
is possible in the case that non-vanishing
Floer cohomology occurs when for each $i=1,\cdots,l$
 the holonomies along discs $D_j$ are equal for
all $j\in I_i$ so that we can define $Q_j$.
 This includes all examples we show in the last section.

\section{A direct computation of ring structure}\label{sec:compute}
Now, we provide two different computations of Floer cohomology rings of torus fibers in toric Fano manifolds.
In this section, we give a direct computation using the classification of holomorphic discs by
Oh and the author in \cite{CO}. For simplicity, we carry out calculations for degree 1 generators, which
is enough to see the whole algebraic structure of the ring due to associativity.

First we choose the generators $C_i$ of $H^1(L)$ for $i=1,\cdots,n$.
\begin{definition}\label{generator}
Let $l_i$ be a cicle $1\times \cdots S^1 \cdots \times 1$ where
$S^1$ is the $i$-th circle of $(S^1)^n \subset (\CC^*)^n $. Then torus action of $(S^1)^n$ on $L$ gives
a corresponding cycles in $L$, which we also denote as $l_i$ by abuse of notation.
For $i=1,\cdots n$, denote by $C_i \in H^1 (L)$ the Poincare dual of the cycle
$$ (-1)^{i-1}(l_1 \times \cdots \times \hat{l_i} \times \cdots l_n).$$
Similarly, We denote by $C_{i,j} \in H^2(L)$ the Poincare dual of
the cycle
$$ (l_1 \times \cdots \times \hat{l_i} \times \cdots \times
\hat{l_j}\times \cdots l_n)$$ for $i \neq j$. and we also define
$C_{i_1,\cdots,i_k}\in H^k(L)$ similarly for the index set
$\{i_1,i_2,\cdots,i_k\}$.
\end{definition}

Now we show that $C_i$'s generate the Floer cohomology ring
$HF^{BM}(L;J_0)$.
\begin{prop}
Let $L$ be a Lagrangian torus fiber whose Floer cohomology group
$HF^{BM}(L;J_0)$ is nonvanishing, thus isomorphic to $H^*(L;\Lambda_{nov})$.
 Then, for each $i$, $C_i$ is a Floer-cycle without any correction
terms, and Floer cohomology $HF^{BM}(L;J_0)$ is generated by $C_i$ for
$i=1,\cdots,n$ as a ring.
\end{prop}
\begin{proof}
From the construction in Definition \ref{correction},
any correction term added to $C_i$, like $PD(C_i) \times Q_j$,
is supposed to have chain dimension $n+1$ or higher. 
Hence, as a current in $L$, they are zero. Hence, $C_i$ itself
is a Floer-cycle.

To see that $\{C_i\}$ generate the Floer cohomology ring, note that
$$m_2(C_{i_1}, m_2(C_{i_2},\cdots, m_2(C_{i_{k-1}}
,C_{i_k})\cdots)$$ is a Floer cycle whose index zero part is
$$m_{2,0}(C_{i_1}, m_{2,0}(C_{i_2},\cdots, m_{2,0}(C_{i_{k-1}}
,C_{i_k})\cdots).$$ Since $m_{2,0}$ is nothing but the cup product,
hence the latter equals $C_{i_1,\cdots,i_k} \in H^k(L)$ up to
sign. Note that all the other terms (terms containing $m_{2,\beta}$ with
non-zero $\beta$) are higher order terms with respect to the filtration by $T$.
Hence, these elements generate the ring $HF^{BM}(L;J_0)$.
\end{proof}
\begin{remark}
For the sign convention for the cup product, see \cite{FOOO} convention 25.14.
\end{remark}

Now, we compute the quantum contribution. We first state the
following lemma which is a special case of Proposition \ref{dim}.
\begin{lemma}\label{twodisc}
Let $\beta \in \pi_2(M,L)$ be a homotopy class.
Then the degree (as a cochain) of $m_{2,\beta} (C_i,C_j)$
is given by $$deg(C_i) + deg(C_j) - \mu(\beta) = 2 - \mu(\beta).$$
Hence, we have non-trivial $m_{2,\beta}$ product between the generators
$C_i$ for $\beta$ with $\mu(\beta)=0$ or $2$.
\end{lemma}
The product when $\mu(\beta)=0$ is the classical cup product, hence
we  consider the contributions from homotopy classes with Maslov index two.
Let us recall the definition of $m_{2,\beta}$.
\begin{equation}
 m_{2,\beta}(C_i,C_j)= (-1)^{n+1}((\MM_3^{\textrm{main}}(\beta_k) \,_{ev_1,ev_2} \times (C_i \times C_j), ev_0),
\end{equation}
where $i$ is an embedding of cycles into $L$.
Recall that {\it main} component is one of the component of moduil space of
discs with marked points $ev_0,ev_1,ev_2$ lie on the disc counter-clockwise
direction. The fact that we use only the main component of the moduli space 
is important, and this make computation a little cumbersome.

To get an intuitive idea about calculations, we first
study the case of $\CP^1$.

\subsection{Example : the equator $L \subset \CP^1$.}
Let $L$ be the equator of $\CP^1$, whose Floer cohomology $HF^(L,L)$ is isomorphic to $H^*(S^1)$.
We pick a point $p$ which will be an element of both the singular
homology $H_0(L)$ and $HF^1(L,J_0)$.
Note that the cup product $$PD(p) \cup PD(p) =0,$$
since generically two points does not intersect in $S^1$.
In our case, we choose $t\in S^1$ which is not equal to $1$, and consider 
two points $p$ and $q = t\cdot p $. 
Then, clearly $$m_{2,0}(p,q) =0.$$

Now, we consider products $m_{2,\beta}$ with non-zero $\beta$.
By Lemma \ref{twodisc}, we only consider $\beta$ with $\mu(\beta)=2$.
Recall from \cite{C} that there are only two such holomorphic discs $D_u,D_l$
(up to $Aut(D^2)$) with boundary on $L$, which are nothing
but discs covering upper(lower)-hemisphere $D_u$ ($D_l$).

\begin{figure}
\begin{center}
\includegraphics[height=2in]{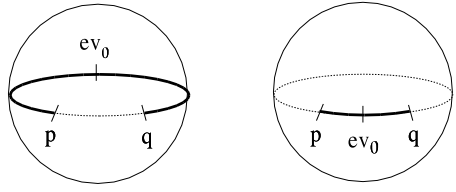}
\caption{Main components of evaluation maps for $D_u$(left) and $D_l$(right)}
\label{fig2}
\end{center}
\end{figure}

Note that both discs intersect $p$ and $q$.
Then, the product $m_{2,D_u} (p,q) $ 
is a certain part of the boundary of $D_u$.
More precisely, since we only consider the ``main'' component 
of the boundary, where $ev_0,p,q$ is ordered counter-clockwise on
the boundary of the disc $D_u$, we obtain a part of $S^1$ as in Figure \ref{fig2}.
And similarly, the product $m_{2,D_l} (p,q) $ only takes the ``main'' component
of the boundary, where $ev_0,p,q$ is ordered counter-clockwise on
the boundary of the disc $D_l$.
Therefore, after adding these two pieces,
we obtain the whole equator:
$$m_2(p,q) = (m_{2,D_u}(p,q) + m_{2,D_l}(p,q))T^{\omega(D)}q = 
[S^1] T^{\omega(D)}q.$$

Hence, $HF(L,J_0) $ is a Clifford algebra with
a generator $[p]$ and a unit $1=[S^1]$ such that,
$$\widetilde{m_2}([p],[p]) = [\widetilde{m_2}(p,q)]= [m_2(p,q)] $$
$$ =[S^1] T^{\omega(D)}q  = 1\cdot T^{\omega(D)}q \in HF(L,J_0)$$

\subsection{Computation of  $m_2(C_i,C_j) + m_2(C_j,C_i)$}
The previous example illustrates that the product $m_2(C_i,C_j)$ is a sum of chains in $L$. Actually
this is not a cycle of singular homology in general as seen in the Definition \ref{correction}.
But we can make the calculation much easier by computing the sum
$m_2(C_i,C_j) + m_2(C_j,C_i)$ instead of each individual pieces.
(Next section generalizes this observation).

The reason that we compute this sum rather than each part, is that
by adding each ``main'' component we will obtain the whole
boundaries of the discs which intersect both $C_i$ and $C_j$.
And computing this sum will be enough to show that the algebra we obtain
is a Clifford Algebra.
Again the only nontrivial $m_{2,\beta}(C_i,C_j)$ will come from 
the homotopy classes $\beta_1,\cdots,\beta_N$ of Maslov index 2
by Lemma \ref{twodisc}.

We recall the relevant fiber product orientation from \cite{C}.
This in fact provides the same orientation as in \cite{FOOO}, which is
described by the orientation of fiber products of Kuranishi structures.
(Smooth simplex may be considered as a weakly submersive strongly
continuous map from a space with Kuranishi structure with
corners where obstruction bundle is taken to be the normal bundle
of the embedding). 
\begin{definition}\cite{C}
Let $X,P,Y$ be an oriented smooth manifolds.
Let $f:X \to Y$ and $i:P \to Y$ be a smooth map.
We define the orientation of the fibre product
$X \,\times_Y P$ for the case that the map $i:P \to Y$
is an embedding.
Let $f : X \to L$ be a submersion and $i : P \to L$ be an embedding.
Here we will regard $P$ as a submanifold of $L$.
By $x, l, p$ we denote the dimension of $X, L, P$.
Take a point $q \in f(X) \cap P$.
We can choose an oriented basis $< u_1, \dots, u_l> \in T_q L$ and 
$< w_1, \dots, w_p> \in T_q P$ which agrees with the given
orientations of $L$ and $P$.
Since $f$ is a submersion, we can choose $< v_1, \dots, v_l> \in T_p X$ for
some $p \in f^{-1} (q)$ such that $(df)_p (v_k) = u_k$ for $k = 1,\dots,l$.
Then, we can choose a basis $< \eta_1, \dots, \eta_{x-l}> \in Ker (df_p)$
such that $< \eta_1, \dots, \eta_{x-l}, v_1, \dots, v_l, > $ is the
given orientation of $T_p X$.
Then we define an orientation on the fibre product $X \,_f \times _i P$ so
that $<\eta_1, \dots, \eta_{x-l},w_1,\dots,w_p>$ becomes an oriented basis.
\end{definition}

From now on, $[\;]$ means the oriented frame on its tangent bundle.
We remark that we mainly follow the amazing work of orientation convention 
in  \cite{FOOO}. 
We may rewrite following \cite{FOOO},
\begin{eqnarray*}
 m_{2,\beta_k}(C_i,C_j) + m_{2,\beta_k}(C_j,C_i)
&=& (-1)^{n+1}((\MM_3(\beta_k) \,_{ev_1,ev_2} \times (C_i \times C_j), ev_0)\\
&=& (-1) ((\MM_3(\beta_k) \,_{ev_1}
\times_{i} C_i)\,_{ev_2} \times_{i} C_j), ev_0)
\end{eqnarray*}
Recall that
\begin{eqnarray*}
[\MM_3(\beta_k)] &=& ([\widetilde{\MM}(\beta_k)]\times[\partial D_0]
\times [\partial D_1] \times [\partial D_2]) / PSL(2:\RR) \\
&=& (-1)([\partial D_0] \times  [\partial D_2]\times
[\widetilde{\MM}(\beta_k)]\times[\partial D_1])/PSL(2:\RR)\\
&=& (-1)([\partial D_0] \times  [\partial D_2]\times[T^n])
\end{eqnarray*}
Here last equality follows from \cite{C} Proposition 3.18.
By the above definition of fiber product orientation, we have
 $$[\MM_3(\beta_k) \,_{ev_1}\times C_i] =
(-1)[\partial D_0] \times  [\partial D_2] \times [C_i]$$
Therefore,
$$ m_{2,\beta_k}(C_i,C_j) + m_{2,\beta_k}(C_j,C_i)
=([\partial D_0] \times  [\partial D_2] \times
[C_i])\,_{ev_2}\times_i [C_j])$$

As the marked point travel around the $3$rd marked point $\partial D_2$,
its trajectory in $L$ is $$v_{k1}l_1 + \cdots+v_{kn} l_n.$$
Here, $v_k$ for $k=1,\cdots,N$ are normal vectors to the
codimension 1 faces of the moment polytope for $M$.
\begin{eqnarray*}
[\partial D_2] \times [C_i] &=& [v_{k1}l_1 +\cdots+v_{kn} l_n ]
 \times (-1)^{i-1}[l_1\times \cdots \times \hat{l_i} \times \cdots \times l_n]
 \\
&=& (-1)^{i-1} [v_{ki} l_i \times
l_1\times \cdots \times \hat{l_i} \times \cdots \times l_n]\\
&=& v_{ki}[l_1 \times \cdots \times l_n] = v_{ki} [T^n]
\end{eqnarray*}
Therefore,
\begin{eqnarray*}
 m_{2,\beta_k}(C_i,C_j) + m_{2,\beta_k} (C_j,C_i) &=&
([\partial D_0] \times  [\partial D_2]
\times [C_i])\,_{ev_2}\times_i [C_j])\\
&=&([\partial D_0]\times [v_{ki}T^n])\,_{ev_2}\times_i [C_j] \\
&=&v_{ki} ([\partial D_0]\times [C_j]) \\
&=& v_{ki}v_{kj} [T^n] \\
&=& v_{ki}v_{kj} \cdot 1
\end{eqnarray*}

Also note that signs of the following the cup product works as
$$m_{2,0} (C_i, C_j) = -m_{2,0}( C_j , C_i).$$ 
Therefore,
\begin{prop}\label{Q}
\begin{eqnarray*}
m_2(C_i,C_j) + m_2 (C_j,C_i)
 &=& \sum_{k=1}^N (m_{2,\beta_k}(C_i,C_j)+  m_{2,\beta_k} (C_j,C_i)) T^{e_k}q \\
 &=& \sum_{k=1}^N v_{ki}v_{kj} T^{e_k}q
\end{eqnarray*}
\end{prop}

Now we consider the case when $i=j$.
The above formula also works for the case $i=j$ 
after we perturb $C_i$ by a torus action to $t \cdot C_i$ for 
some $t \in T^n$.
Also we have the following easy lemma.
\begin{lemma}
$$[m_2(C_i,t \cdot C_i)] = [m_2(t \cdot C_i,C_i)]\;\; \textrm{in}\; HF^*(L,J_0),$$
\end{lemma}

\begin{corollary}\label{pro2}
 $$m_2(C_i,t \cdot C_i) = \sum_{k=1}^N \frac{1}{2} v_{ki}^2 
\otimes T^{e_k}q.$$
\end{corollary}

Now, we recall the definition of the Clifford algebra.
\begin{definition}
Let $V$ be a $\QQ$-vector space with
a non-degenerate symmetric bilinear form $Q$ on $V$.
The Clifford Algebra $Cl(V,Q)$ is defined as
$$ Cl(V,Q) = T(V)/I(V,Q),$$
where $T(V)$ is the tensor algebra
$$T(V) = \bigoplus_{k=0} V^k,$$
And $I(V,Q)$ is the ideal in $T(V)$ generated by elements
$$ v \otimes v -\frac{1}{2} Q(v,v)1 \;\; \textrm{for}\;\; v \in V.$$
Alternatively, one may define Cl(V,Q) with the
relation
$$v \cdot w + w \cdot v = Q(v,w).$$
\end{definition}

In our case, we consider a universal Novikov ring $\Lambda_{nov}$ instead of
$\QQ$ as a coefficient.
Now Proposition \ref{Q} and Corollary \ref{pro2} implies our main theorem.

\begin{theorem}\label{main}
Let $L \subset M$ be a Lagrangian torus fiber in Fano toric manifold whose Floer cohomology is non-vanishing.
Then, the Floer cohomology ring $(HF^{BM}(L;J_0), \widetilde{m_2})$
has a Clifford Algebra structure with
generators given by $C_i$ for $i=1,\cdots, n$, and
its relations as 
$$ \widetilde{m_2}( C_i,C_j) +  \widetilde{m_2}( C_j,C_i) = Q(C_i,C_j),$$
where symmetric bilinear form $Q$ is given by 
$$ Q(C_i,C_j) = 
 \sum_{k=1}^N v_{ki}v_{kj} T^{e_k}q$$
Furthermore, this $Q$ agrees with the Hessian of the superpotential 
$W(\Theta)$ of the mirror Landau-Ginzburg model of toric Fano manifold.
(upon the substitution ``$T^{2\pi} = e^{-1}$'').
\end{theorem}
\begin{proof}
We only need to check the last statement.
Recall that the superpotential  is given as (see for example \cite{HV},
\cite{CO})
$$W(\Theta) = \sum_{k=1}^N e^{-y_k - <\Theta,v_k>}$$
Hence, it is easy to see that
$$\frac{\partial W(\Theta)}{\partial \Theta_i} = -
\sum_{k=1}^N v_{ki} e^{-y_k - <\Theta,v_k>}  $$
And 
$$\frac{\partial^2 W(\Theta)}{\partial \Theta_i \partial \Theta_j} =
\sum_{k=1}^N v_{ki} v_{kj} e^{-y_k - <\Theta,v_k>}  $$
Here $\Theta$ is a coordinate on mirror Landau-Ginzburg model, and it is
related to the toric manifold $M$ as follows.
Real part of the variable $\Theta$ is given by $(a_1,\cdots,a_n) \in P$ which 
is the image point of the Lagrangian torus fiber $L$ in moment polytope $P$,
whereas imaginary part is given by holonomy of the flat line bundle along $L$.
When the Floer cohomology of $L$ is non-vanishing, 
the corresponding $\Theta$ becomes the critical point of $W$ as shown 
in \cite{CO}, and $2\pi$ times its exponent $(y_k + <\Theta,v_k>)$ becomes 
the area of holomorphic discs which we denoted as $e_k$ in this paper.
If we ignore harmless grading $q$, and with the equivalence ``$T^{2\pi} = e^{-1}$'',
$$ e^{-y_k - <\Theta,v_k>} = T^{e_k}.$$
Hence, this proves the claim.
\end{proof}
More correspondences will be given in the next section.

\section{Analogue of divisor equation for discs.}\label{divisor}
In this section, we introduce an analogue of divisor equation and this will explain how Clifford algebra structure
naturally arises for Floer cohomology rings of Lagrangian submanifolds, as this section provides the alternative proof
of results in the previous section. To state the result, it is better to write down the
formula in terms of $\LI$-algebra (strong homotopy Lie algebra) maps. Recall that every $\AI$-algebra has an underlying $\LI$-algebra structure by the  following relation (This is similar to the fact that commutator of an associative algebra $A$
defines a Lie algebra on $A$). 

\begin{theorem}[\cite{LM},\cite{LS}](or see \cite{Fuk2}).
An $\AI$-structure $\{ m_k: \otimes^k V \to V \}$ on the graded vector space $V$ induces an
$\LI$-structure $\{ l_k: \otimes^k V \to V \}$ where for all non-negative integer $k$, $\beta \in \pi_2(M,L)$, 
\begin{equation}\label{li}
 l_{k,\beta}( v_1 \otimes \cdots \otimes v_k ) = \sum_{\sigma \in S_n} (-1)^{\epsilon(\sigma)}
m_{k,\beta}(v_{\sigma(1)} \otimes \cdots \otimes v_{\sigma(k)}),  
\end{equation}
$$\textrm{with} \; \epsilon(\sigma) = \sum_{i,j \;\textrm{with}\; i<j, \sigma(i) > \sigma(j)} (deg(v_i) +1 )(deg(v_j)+1).$$
\end{theorem}
Namely, $l_k$ map is a skew-symmetrization of $m_k$ map. 
For example, $$l_{2,\beta}(x,y) = m_{2,\beta}(x,y) + (-1)^{(x+1)(y+1)} m_{2,\beta}(y,x). $$

The following is the Divisor equation of Gromov-Witten invariants. For a general equation involving
gravitational descendents, see \cite{H0}.
\begin{prop}[\cite{KM}]
Let $M$ be a convex algebraic manifold, For $\alpha \in H_2(M)$, let $I_{g,m,\alpha}^M : H^*(V)^{\otimes n} \to
H^{*}(\overline{M}_{g,n})$ be the Gromov-Witten invairants. 
Then, for $\gamma_1 \in H^2(M)$, and $\pi_n:\overline{M}_{g,n} \to \overline{M}_{g,n-1}$, we have
$$\pi_{n*}(I_{g,n,\alpha}^M(\gamma_1 \otimes \cdots \otimes \gamma_n) = (\alpha \cdot \gamma_1)
I_{g,n-1,\alpha}^M(\gamma_2 \otimes \cdots \otimes \gamma_n).$$
\end{prop}

 Now, we state an {\bf analogue of  divisor equation for discs}. 

\begin{prop}\label{prop:div}
If $P_i$ is a cycle of cohomology degree 1 in $L$, then for $k \geq 1$,
\begin{equation}
l_{k,\beta}(P_1,\cdots,P_k)  = (P_i \cdot \partial \beta) \; l_{k-1,\beta}(P_1,\cdots,\widehat{P_i},\cdots,P_k),
\end{equation}
where $\partial :\pi_2(M,L) \to \pi_1(L)$, and $\widehat{P_{i}}$ means that $P_i$ term is omitted. 
\end{prop}
\begin{remark}
Here is the sign convention for intersection of two chains $P,Q$ of complementary degree in $L$. 
At each transversal intersection $p \in P \cap Q$, for a basis $[T_pP]$ of tangent space $T_p P$, 
and similarly for $[T_p Q]$ and $[T_p L]$, if $[T_pP][T_pQ]$ has the same orientation as $[T_pL]$ then
it is counted as $(+1)$, otherwise it is counted as $(-1)$.
\end{remark}
Before we prove the proposition, we show how to prove the results in the previous section using the analogue of divisor equation for discs.
Recall that by $\beta_k \in \pi_2(M,L)$ for $k=1,\cdots,N$, we denote the homotopy class of holomorphic disc of Maslov index two corresponding to $N$ codimension one facets of the moment polytope(\cite{CO}).

By definition, we have
\begin{eqnarray*}
l_{0,\beta_k} &=& m_{0,\beta_k} = T^{e_k}q \\
l_{1,\beta_k}(P) &=& m_{1,\beta_k}(P)
\end{eqnarray*}

For a degree 1 generators $C_i, C_j$ of $H^*(L)$ which defined in Definition \ref{generator}, 
we apply the divisor equation for discs repeatedly
\begin{eqnarray*}
l_{2,\beta_k}(C_i,C_j) &=& (C_i \cdot \partial \beta_k) l_{1,\beta_k}(C_j) \\
&=& (C_i \cdot \partial \beta_k)(C_j \cdot \partial \beta_k) l_{0,\beta_k} \\
&=& (v_{ki})(v_{kj})\otimes T^{e_k}q
\end{eqnarray*}
The last equality follows from the definitions that 
$$C_i = (-1)^{i-1}(l_1 \times \cdots \times \hat{l_i} \times \cdots l_n),$$
$$\partial \beta_k = v_{k1}l_1 + \cdots + v_{kn}l_n.$$
Hence, it is easy to see that 
$$C_i \cdot \partial \beta_k = (-1)^{n}v_{ki}.$$
Hence, we obtain the lemma \ref{Q}, as we have $l_2(C_i,C_j) = m_2(C_i,C_j) + m_2(C_j,C_i)$. 

In general, we have
\begin{corollary} For any Lagrangian torus fiber $L$ in toric Fano manifold $M$(
 whose Floer cohomology may be vanishing), we have
\begin{eqnarray*}
l_m(C_{i_1},\cdots,C_{i_m}) &=&
\sum_{k=1}^N l_{m,\beta_k}(C_{i_1},\cdots,C_{i_m}) \\
&=& (-1)^{nm}\sum_{k=1}^N v_{ki_1} \cdots v_{ki_m}  \otimes T^{e_k}q \\
&=& (-1)^{(n-1)m}\frac{\partial^m W(\Theta)}{\partial \Theta_{i_1} \cdots \partial \Theta_{i_m}},
\end{eqnarray*}
where $W(\Theta)$ is the superpotential of Landau-Ginzburg mirror model of $M$. 
\end{corollary}
This corollary extends the correspondence observed in \cite{CO}, $m_0 = l_0 = W(\Theta)$. 
Note that such correspondence considered at every Lagragian torus fibers with flat line bundles may
be used to recover the superpotential $W(\Theta)$ of the Landau-Ginzburg mirror. But the above corollary indicates that in fact
one Lagrangian torus fiber with a fixed flat line bundle (whose Floer cohomology may
be vanishing) in $M$ is enough to recover the superpotential in this case: It is because  the superpotential is a holomorphic function on $(\CC^*)^n$ and all its partial derivatives at the corresponding point on the mirror is
given from the products of $\LI$-algebra by the above correspondence.

Also note that above product does not depend on the choice of cycles $C_*$ 
since it is determined by the intersection numbers which only depends on the homology class of $C_*$. 
These are also invariants with respect to the change of an almost complex structure.
By $J_0$ we denote the standard complex structure of toric Fano manifold $M$,
and denote the corresponding $l_m$ products by $l^{J_0}_m$. 
\begin{prop}\label{invariant}
Let $L$ be any Lagrangian torus fiber of toric Fano manifold $M$. 
Let $J_1 \in J_{reg}(M)$ be a tame almost complex structure such that all simple $J$-holomorphic 
discs are Fredholm regular. 
Then, for $k =1 ,\cdots, N$, we have 
$$l^{J_0}_{m,\beta_k}(C_{i_1},\cdots,C_{i_m}) = l^{J_1}_{m,\beta_k}(C_{i_1},\cdots,C_{i_m})$$ in $H^*(L;\Lambda_{nov}).$
\end{prop}
\begin{proof}
As in \cite{MS}, one can prove that the subset $J_{reg}(M)$ is of second category and path connected.
Since any $J$-holomorphic disc with Maslov index two is simple and its homotopy class $\beta_k$ is minimal, the moduli space $\CM(\beta_k;J_t)$ of $J_t$ holomorphic discs is in fact a manifold without boundary.
Then, by choosing a path $J_{t} \in J_{reg}(M)$,
we set $$\CM_{m+1}(\beta_k;\mathcal{J}) = \cup_{t \in [0,1]} \big( \{t\} \times \CM_{m+1}(\beta_k;J_t) \big)$$
Then we have 
$$\partial \big( \CM_{m+1} (\beta_k;\mathcal{J})_{ev} \times (\prod_{j=1}^m C_{i_m})\big) = 
\CM(\beta_k;J_1)_{ev} \times (\prod_{j=1}^m C_{i_m}) - \CM(\beta_k;J_0)_{ev} \times (\prod_{j=1}^m C_{i_m}) $$
which proves the proposition.
\end{proof}
Now we begin the proof of the proposition \ref{prop:div}
\begin{proof}
Rough idea is that if $P_i$ is a cycle of codimension 1, then it always intersects with the boundary
of a $J$-holomorphic disc of homotopy class $\beta$
with $(P_i \cdot \partial \beta)$ number of times (counted with sign). Hence, if $P_i$ is dropped from
the argument of $m_k$, the resulting image should be the same up to a multiple of the intersection number.
While this is the same idea as the ``divisor equation '' in Gromov-Witten theory, there are a few differences.
First, $m_k$ map records only part of the boundaries of J-holomorphic discs
 as it is defined by using only the main component $\CM_{k}^{main}$. But note that intersection
of $P_i$ and the disc may occur at arbitrary point of the domain $\partial D^2$. 
Hence we consider $\LI$-algebra map, $l_k$, which will be shown to record
the whole boundaries of discs. Then, the next step involves delicate sign
analysis in the case that the parameter $P_i$ is dropped from the $m_k(P_1,\cdots, P_k)$ to
obtain $m_{k-1}(P_1,\cdots,\widehat{P_i},\cdots,P_k)$ for a codimension 1 cycle $P_i$.

Suppose there exist an element $((D^2 ,\vec{z}),h) \in \CM_{k+1}^{main}(\beta)$, where
$h:D^2 \to M$ is a $J$-holomorphic map. We also assume that for a fixed chains
$P_1,\cdots, P_k$ in $L$, we have $h(z_i) \in P_i$ for each $i = 1,\dots,k$. 
Boundary marked points $z_1,\cdots,z_k$ ($z_0$ is omitted here)
seperates $\partial D^2$ into $k$ connected pieces. And only the component between
$k$-th and $1$-th marked point contributes to the chain $m_k(P_1,\cdots,P_k)$,
as it is obtained as an evaluation of $0$-th marked point which lies between 
those two marked point in $\CM_{k+1}^{main}$.
Now, it is easy to see that up to sign, other connected components 
will contribute to the chains 
$$m_k(P_2,\cdots,P_k,P_1), m_k(P_3,\cdots,P_1,P_2), \cdots, m_k(P_k, P_1, \cdots,P_{k-1}).$$
and $((D^2 ,\vec{z}),h)$ will not contribute to other terms of $l_k(P_1,\cdots,P_k)$ generically
due to the ordering of marked points. Now, we show that signs in (\ref{li}) is needed 
to have a coherent sign in the images of the above chains.

We recall the following lemma from \cite{FOOO}.
\begin{lemma}[FOOO, Lemma 25.3]\label{sign:f5}
Let $\sigma$ be the transposition element $(i,i+1)$ in the $k$-th
symmetric group $S_k$. Then the action of $\sigma$ on 
$\CM_1(\beta,P_1,\cdots,P_i,P_{i+1},\cdots,P_k)$
by changing the order of marked points is described by the
following.
\begin{equation*}
\sigma(\CM_1(\beta,P_1,\cdots,P_i,P_{i+1},\cdots,P_k))
\end{equation*}
$$ = (-1)^{(deg\;P_i +1)(deg\;P_{i+1} +1)}
\CM^{\sigma}_1(\beta,P_1,\cdots,P_{i+1},P_i,\cdots,P_k)$$
\end{lemma}
\begin{remark}
In the first term, $\CM_1(\beta,P_1,\cdots,P_i,P_{i+1},\cdots,P_k))$ is defined by using the moduli
space with boundary marked points lying cyclically, whereas in the second term,
$\CM^\sigma_1(\beta,P_1,\cdots,P_{i+1},P_i,\cdots,P_k)$ is defined by using the moduli space $\CM_k^{\sigma}$
with boundary marked points lying in the order $z_0,\cdots,z_{i-1},z_{i+1},z_{i},z_{i+2},\cdots,z_k$.
Namely, in the latter case, only the labeling of two marked point is changed from the first case.
\end{remark}
Let $\sigma \in S_n$ be a permutation denoted by $(1,2,\cdots, k)$. (i.e. $1 \to 2, 2\to 3, \cdots, k \to 1$). Then, by applying the above lemma repeatedly, we have
$$\sigma(\CM_1(\beta,P_1,\cdots,P_k)) = (-1)^{\epsilon(\sigma)} \CM^\sigma_1(\beta,P_2,\cdots,P_k,P_1),$$
where $\epsilon(\sigma)$ is the same sign as appeared in (\ref{li}).
Now, it is not hard to check that  the latter has the same sign as $(-1)^{\epsilon(\sigma)}\CM(\beta;P_2,\cdots,P_k,P_1)$.
Here, $\CM^\sigma_1(\beta,P_2,\cdots,P_k,P_1)$ and $\CM(\beta;P_2,\cdots,P_k,P_1)$ have different images comming
from the same set of $J$-holomorphic discs. This is because that marked points in the former case lie on the circle 
in the order $0,k,1,2,\cdots,k-1$ and in the latter case marked points lie on the circle
in the order $k,0,1,2,\cdots,k-1$. Hence, as we evaluate at 0-th marked point, their images come
from the neighboring connected components of $\partial D^2$ seperated by marked points.
Hence, this proves that with the sign given as in (\ref{li}), the image of the boundary of discs can be glued 
in the $l_k$ map.

Now, we explain the second step which computes the change of sign
as the argument $P_i$ is dropped from the $m_k(P_1,\cdots, P_k)$ to
obtain $m_{k-1}(P_1,\cdots,\widehat{P_i},\cdots,P_k)$ for a codimension 1 cycle $P_i$.
In the computation, we will calculate the $i$-th fiber product with $P_i$ to
remove the term from the fiber product. 

Let $((D^2 ,\vec{z}),h) \in \CM_{k}^{main}(\beta)$, where
$h:D^2 \to M$ is a $J$-holomorphic map. We also assume that for a fixed chains
$P_1,\cdots,\widehat{P_i},\cdots, P_k$ in $L$, we have $h(z_i) \in P_i$ for each $i = 1,\dots,\widehat{i},\dots,k$.
And let $P_i$ be a cycle of codimension one in $L$. If $[P_i] \cdot \partial \beta$ is not zero, 
then, a generic cycle $P_i$ should intersect with $ h(\partial D^2)$ transversally.
Hence, we obtain a corresponding element $((D^2 ,\vec{z'}),h) \in \CM_{k+1}^{main}(\beta)$
with $h(z_i') \in P_i$ for each $i = 1,\dots,k$.

We recall that the moduli space $\CM_{k+1}^{main}(\beta)$ is oriented as
$$ \big( [\widetilde{\CM}(\beta)] \times [\partial D^2_0] \times \cdots \times [\partial D_k^2] \big) / PSL(2;\CC), $$
where $[\partial D^2_i]$ denotes the tangent vector corresponding to the counterclockwise rotation of $i$-th marked point.
If we take $[\partial D^2_i]$ to the last, we have
$$ = (-1)^{s_1} \big(\big( [\widetilde{\CM}(\beta)] \times [\partial D^2_0] \times \cdots 
\widehat{[\partial D^2_i]} \times [\partial D_k^2] \big) / PSL(2;\CC) \big)\times [\partial D^2_i]$$
where $s_1 = k-i +1$.
Now we write
$$ \CM(\beta;P_1,\cdots,P_k)
= (-1)^{s_2} \CM_{k+1}^{main} (\beta)_{ev_1,\dots,ev_k} \times (P_1 \times \cdots \times P_k)$$
$$= (-1)^{s_3} \big( \cdots \big( \CM_{k+1}^{main} (\beta)_{ev_1} \times P_1\big) \cdots_{ev_{k}} \times P_k \big),$$
where $s_2 = (n+1)\sum_{l=1}^{k-1} \sum_{j=1}^l deg (P_j), s_3 = \sum_{l=1}^{k-1} \sum_{j=1}^l deg (P_j)$.

Then, if we look at the term
$\big( \CM_{k+1}^{main} (\beta)_{ev_1} \times P_1\big)$, it can be oriented as
$$ (-1)^{s_1}\big(\big(\big( [\widetilde{\CM}(\beta)] \times [\partial D^2_0] \times \cdots 
\widehat{[\partial D^2_i]} \times \cdots \times [\partial D_k^2] \big) / PSL(2;\CC) \big)\times [\partial D^2_i] \big)_{ev_1} \times [P_1]$$
$$= (-1)^{s_4} \big(\big(\big( [\widetilde{\CM}^o(\beta)] \times [\partial D^2_0] \times \cdots 
\widehat{[\partial D^2_i]} \times \cdots \times [\partial D_k^2] \big) / PSL(2;\CC) \big)\times [\partial D^2_i] \times [L] \big)_{ev_1} \times [P_1]$$
$$ = (-1)^{s_4} (\big(\big( [\widetilde{\CM}^o(\beta)] \times [\partial D^2_0] \times \cdots 
\widehat{[\partial D^2_i]} \times \cdots \times [\partial D_k^2] \big) / PSL(2;\CC) \big)\times [\partial D^2_i] \times [P_1]$$
$$ = (-1)^{s_5}(\big(\big( [\widetilde{\CM}^o(\beta)] \times [\partial D^2_0] \times \cdots 
\widehat{[\partial D^2_i]} \times \cdots \times[\partial D_k^2] \big) / PSL(2;\CC) \big) \times [P_1] \times [\partial D^2_i]$$
$$ = (-1)^{s_6} 
 \big(\big(\big( [\widetilde{\CM}(\beta)] \times [\partial D^2_0] \times \cdots 
\widehat{[\partial D^2_i]} \times \cdots \times [\partial D_k^2] \big) / PSL(2;\CC) \big)_{ev_1} \times [P_1] \big) \times [\partial D^2_i]$$
where $[\widetilde{\CM}^o(\beta)][L] = [\widetilde{\CM}(\beta)]$ and 
$s_4 = n(k+2) + s_1 $, $s_5 =  s_4 + p_1 = s_4 + dim(P_1)$ and $s_6 = s_5 + n(k+1)$. 
Hence $s_6 = n + p_1 + k -i +1  = deg (p_1) + k - i +1$
Now we repeat this process up to $P_{i-1}$ and 
$$\big( \cdots \big( \CM_{k+1}^{main} (\beta)_{ev_1} \times P_1\big) \cdots \times P_{i-1} \big)$$
is oriented as
$$ (-1)^{s_7} 
 \big(\cdots (\big( [\widetilde{\CM}(\beta)] \times [\partial D^2_0] \times \cdots 
\widehat{[\partial D^2_i]} \times \cdots \times [\partial D_k^2] \big) / [PSL(2;\CC)] \big)_{ev_1} $$
$$ \times [P_1] \big) \times \cdots \times [P_{i-1}] \big) \times [\partial D^2_i]$$
with $s_7 = deg (P_2) + \cdots + deg (P_{i-1}).$ Then,
$$\big( \cdots \big( \CM_{k+1}^{main} (\beta)_{ev_1} \times P_1\big) \cdots \times P_{i-1} \big)_{ev_{i}} \times P_i$$
is oriented as
$$ (-1)^{s_8}
 \big(\cdots (\big( [\widetilde{\CM}(\beta)] \times [\partial D^2_0] \times \cdots 
\widehat{[\partial D^2_i]} \times \cdots \times [\partial D_k^2] \big) / [PSL(2;\CC)] \big)_{ev_1} $$
$$ \times [P_1] \big) \times \cdots \times [P_{i-1}] \big)^o \times [\partial D^2_i] [P_i] $$
$$ = (-1)^{s_9}
 \big(\cdots (\big( [\widetilde{\CM}(\beta)] \times [\partial D^2_0] \times \cdots 
\widehat{[\partial D^2_i]} \times \cdots \times [\partial D_k^2] \big) / [PSL(2;\CC)] \big)_{ev_1} $$
$$ \times [P_1] \big) \times \cdots \times [P_{i-1}] \big)^o \times  [P_i] [\partial D^2_i] $$
$$ = (-1)^{s_{10}}
 \big(\cdots (\big( [\widetilde{\CM}(\beta)] \times [\partial D^2_0] \times \cdots 
\widehat{[\partial D^2_i]} \times \cdots \times [\partial D_k^2] \big) / [PSL(2;\CC)] \big)_{ev_1} $$
$$ \times [P_1] \big) \times \cdots \times [P_{i-1}] \big), $$
where 
$s_8 = \sum_{j=1}^{i-1} deg (P_j) + n,\; s_9 = s_8 + (n-1)\cdot 1, \; s_{10} = s_{9} + \epsilon$. The last equality 
follows from the sign of the intersection $[P_i][\partial \beta] = (-1)^{\epsilon} [L]$. 
Hence $s_{10} = \epsilon + (n-1) + n + \sum_{j=1}^{i-1} deg (P_j) + k -i + 1$.

Now, the last expression can be considered as an orientation of 
$$\big( \cdots \big( \CM_{k}^{main} (\beta)_{ev_1} \times P_1\big) \cdots_{ev_{i-1}} \times P_{i-1} \big).$$
Hence, orientation of $$\CM(\beta;P_1,\cdots,P_k)$$ corresponds to
$$ (-1)^{s_3} \big( \cdots \big( \CM_{k+1}^{main} (\beta)_{ev_1} \times P_1\big) \cdots_{ev_{k}} \times P_k \big),$$
$$ \subset (-1)^{s_{11}} \big( \cdots \big( \CM_{k}^{main} (\beta)_{ev_1} \times P_1\big) 
\cdots \times \widehat{P_i} \big) \times \cdots_{ev_{k}} \times P_k \big),$$
$$ = (-1)^{s_{12}} \CM(\beta;P_1,\cdots,\widehat{P_i},\cdots,P_k),$$
where $s_{11} = s_3 + s_{10} $, $s_{12}$ is obtained in a similar way as $s_3$ and
we have $s_{12} = \epsilon + k-i + (k-i)deg (P_i) = \epsilon$, since $deg (P_i) = 1$. 

This proves that if $P_i$ intersect with the $J$-holomorphic disc at several boundary points,
then at each intersection, the sign change  between $m_k(P_1,\cdots, P_k)$ and
$m_{k-1}(P_1,\cdots,\widehat{P_i},\cdots, P_k)$ of the contribution
from this $J$-holomorphic disc, is given by $(-1)^{\epsilon}$ where
$\epsilon$ is the sign of the intersection between $[P_i]$ and the $J$-holomorphic disc at each intersection point. 

Now, we prove the proposition. Note that in the expression $l_k(P_1,\cdots,P_k)$, a term $m_k(P_1,\cdots,P_k)$ carries the same sign as the term $m_k(P_1,\cdots,\widehat{P_i},\cdots,P_k,P_i)$ or any other term which is obtained by moving $P_i$ around if $deg(P_i) =1$.
 Consider the $J$-holomorphic disc contributing non-trivially to the expression $m_{k-1,\beta}(P_1,\cdots,\widehat{P_i},\cdots,P_k)$. Then whole boundary of this disc contribute to $l_{k-1,\beta}(P_1,\cdots,\widehat{P_i},\cdots,P_k)$.
Generically, $P_i$ may intersect with $\partial \beta$ at arbitrary points of the domain $\partial D^2$. 
Suppose such disc intersect $P_i$ between $j$ and $j+1$-th marked point with intersection sign $(-1)^{\epsilon}$ for $j>i+1$
without loss of generality. Then, the whole boundary of this disc would contribute to the terms (while divided into
several pieces)
$$m_{k,\beta}(P_1,\cdots,\widehat{P_i},\cdots, P_{j-1},P_i,P_j,\cdots, P_k), 
+ m_{k,\beta}(P_2,\cdots,\widehat{P_i},\cdots, P_{j-1},P_i,P_j,\cdots, P_k, P_1)$$
$$ + \cdots +  m_{k,\beta}(P_k, P_1, \cdots,\widehat{P_i},\cdots, P_{j-1},P_i,P_j,\cdots, P_{k-1}).$$
By applying the sign analysis, the above terms correspond to terms in $l_k(P_1,\cdots,\widehat{P_i},\cdots,P_k)$
with multiplicity $(-1)^{\epsilon}$. By adding up all the possibilities of intersections between $P_i$ and
$\partial \beta$, we obtain the proposition.
\end{proof}

\section{examples}
In what follows we omit area terms $T^{e_i}q$ for simplicity.
\subsection{The Clifford torus $T^2 \subset \CP^2$}
In \cite{CO}, it is shown that the Clifford torus is the only Lagrangian torus
fiber whose Floer cohomology is non-vanishing, which is isomorphic to
$H^*(T^2;\Lambda_{nov})$.
Hence, by Theorem \ref{main}, $HF^*(T^2,T^2)$ as a ring is a Clifford algebra with
two generators $C_1$ and $C_2$.
Using its moment polytope data, one can immediately compute the matrix of the symmetric bilinear
form 
$$ Q = \left( \begin{array}{cc}
2 & 1/2  \\
1/2 & 2 
\end{array} \right).$$
(See \cite{KL} for computations of $B$-model by physical arguments and the predictions
made for the Clifford torus case.)

\begin{figure}[ht]
\begin{center}
\includegraphics[height=1.5in]{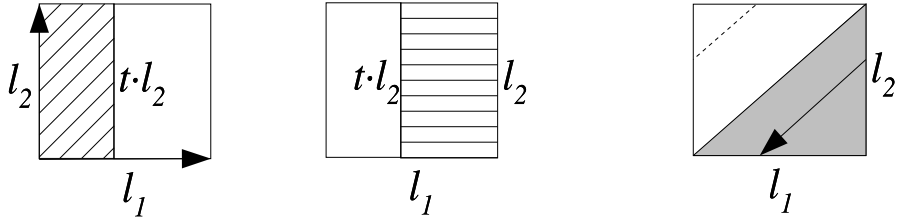}
\caption{$m_{2,\beta_0}(l_2,tl_2) + m_{2,\beta_1}(l_2,tl_2)=L,\;\; m_{2,\beta_0}(l_1,l_2)$}
\label{fig3}
\end{center}
\end{figure}

But it is also instructive to compute $m_2(C_1,C_1)$ and $m_2(C_1,C_2)$ directly.
Consider $T^2$ as a rectangle whose edges are glued accordingly.
Let us assume that its edges are cycles $l_1,l_2$ as given in the Definition
\ref{generator}.
Then by definition, we have $$C_1 = l_2, C_2 = -l_1.$$
First we consider $m_2(C_1,C_1) = m_2(l_2,l_2)$.
As before, we pick $t \in T^2$ so that $l_2$ and $tl_2$ do not intersect.
Then, $$m_{2,0}(l_2,tl_2) =0.$$
Recall that there exists 3 holomorphic discs (up to $Aut(D^2)$) with 
boundary trajectory as 
$$ \partial D_0 = -l_1-l_2,\; \partial D_1 = l_1,\; \partial D_2 = l_2.$$
For $m_{2,\beta}(l_2, tl_2)$, holomorphic discs $D_0$, $D_1$ contributes nontrivially.
Since we only consider {\it main} components as in Figure \ref{fig3}, we have
$$m_{2,\beta_0}(l_2,tl_2) + m_{2,\beta_1}(l_2,tl_2) = [L].$$
Therefore, we have 
$$ \widetilde{m_2}([C_1],[C_1]) = [m_2(C_1,tC_1)] = [L] T^{\omega(D)}q.$$
This agrees with the Corollary \ref{pro2}.
From the Figure \ref{fig3}, it is easy to see that 
the product $m_2(l_2,tl_2)$ is independent of $t\in  T^2$.

Now, we consider the product $m_2(C_1,C_2)$. 
$$m_2(C_1,C_2) = m_{2,0}(C_1,C_2) + \sum_\beta m_{2,\beta}(C_1,C_2) T^{Area(\beta)}q$$
Here $m_{2,0}(C_1,C_2)$ is a cup product which is nothing but
the Poincare dual of intersection $C_1 \cap C_2 = point$.
But as we discussed in Example \ref{pt}, point itself is not
a Floer-cycle. The needed correction term $Q$ 
is obtained in this case from the quantum contribution $m_{2,\beta}$.

It is easy to see that only $\beta_0$ disc contributes to
the product $m_{2,\beta}$, since other discs generically do not
intersect both $C_1,C_2$.
Now, $m_{2,\beta_0}(C_1,C_2)$
 is not cycle but a chain as drawn in Figure \ref{fig3}, since we
only evaluate on the main component.
This is the chain $Q$ that we added to make $<pt>$ a Floer cycle in Definition \ref{correction}.
$$m_2(C_1,C_2) = m_{2,0}(C_1,C_2) + m_{2,\beta_0}(C_1,C_2) T^eq$$
$$ = <pt>  + Q T^e q$$

\subsection{$\CP^1 \times \CP^1$}
Consider $\CP^1 \times \CP^1$ whose moment map image is a rectangle.
For the equator $S^1 \in \CP^1$, 
$S^1 \times S^1 \subset \CP^1 \times \CP^1$ has
nontrivial Floer cohomology and its product
structure is given by 
$$ Q = \left( \begin{array}{cc}
2 & 0  \\
 0& 2
\end{array} \right)$$

\subsection{$\CP^n$} 
The example $\CP^2$ easily generalizes to $\CP^n$. 
The Floer cohomology of the Clifford torus $T^n \subset \CP^n$ 
becomes the Clifford Algebra with $n$ generators with 
symmetric bilinear form as 
$$ Q = \left( \begin{array}{cccc}
2 & 1/2 & \cdots & 1/2  \\
1/2 & 2&  \cdots & 1/2 \\
 \vdots & \vdots & \ddots & \vdots \\
1/2 & 1/2 & \cdots & 2 
\end{array} \right)$$

\bibliographystyle{amsalpha}

\end{document}